\documentclass[12pt]{article} 
\usepackage{latexsym} 
\usepackage{amsmath}
\usepackage{amssymb} 
\usepackage{amscd} 
\usepackage{graphicx}
\usepackage{layout}
\begin{document} 
\newtheorem{Th}{Theorem}[section]
\newtheorem{Cor}{Corollary}[section]
\newtheorem{Prop}{Proposition}[section]
\newtheorem{Lem}{Lemma}[section]
\newtheorem{Def}{Definition}[section]
\newtheorem{Rem}{Remark}[section]
\newtheorem{Ex}{Example}[section]
\newtheorem{stw}{Proposition}[section]


\newcommand{\bet}{\begin{Th}}
\newcommand{\ent}{\stepcounter{Cor}
   \stepcounter{Prop}\stepcounter{Lem}\stepcounter{Def}
   \stepcounter{Rem}\stepcounter{Ex}\end{Th}}


\newcommand{\bec}{\begin{Cor}}
\newcommand{\enc}{\stepcounter{Th}
   \stepcounter{Prop}\stepcounter{Lem}\stepcounter{Def}
   \stepcounter{Rem}\stepcounter{Ex}\end{Cor}}
\newcommand{\bep}{\begin{Prop}}
\newcommand{\enp}{\stepcounter{Th}
   \stepcounter{Cor}\stepcounter{Lem}\stepcounter{Def}
   \stepcounter{Rem}\stepcounter{Ex}\end{Prop}}
\newcommand{\bel}{\begin{Lem}}
\newcommand{\enl}{\stepcounter{Th}
   \stepcounter{Cor}\stepcounter{Prop}\stepcounter{Def}
   \stepcounter{Rem}\stepcounter{Ex}\end{Lem}}
\newcommand{\bef}{\begin{Def}}
\newcommand{\enf}{\stepcounter{Th}
   \stepcounter{Cor}\stepcounter{Prop}\stepcounter{Lem}
   \stepcounter{Rem}\stepcounter{Ex}\end{Def}}
\newcommand{\ber}{\begin{Rem}}
\newcommand{\enr}{
   \stepcounter{Th}\stepcounter{Cor}\stepcounter{Prop}
   \stepcounter{Lem}\stepcounter{Def}\stepcounter{Ex}\end{Rem}}
\newcommand{\bee}{\begin{Ex}}
\newcommand{\ene}{
   \stepcounter{Th}\stepcounter{Cor}\stepcounter{Prop}
   \stepcounter{Lem}\stepcounter{Def}\stepcounter{Rem}\end{Ex}}
\newcommand{\Proof}{\noindent{\it Proof\,}:\ }

\newcommand{\EE}{\mathbf{E}}
\newcommand{\QQ}{\mathbf{Q}}
\newcommand{\R}{\mathbf{R}}
\newcommand{\C}{\mathbf{C}}
\newcommand{\ZZ}{\mathbf{Z}}
\newcommand{\KK}{\mathbf{K}}
\newcommand{\NN}{\mathbf{N}}
\newcommand{\PP}{\mathbf{P}}
\newcommand{\HH}{\mathbf{H}}
\newcommand{\uuu}{\boldsymbol{u}}
\newcommand{\xxx}{\boldsymbol{x}}
\newcommand{\aaa}{\boldsymbol{a}}
\newcommand{\bbb}{\boldsymbol{b}}
\newcommand{\AAA}{\mathbf{A}}
\newcommand{\BBB}{\mathbf{B}}
\newcommand{\ccc}{\boldsymbol{c}}
\newcommand{\iii}{\boldsymbol{i}}
\newcommand{\jjj}{\boldsymbol{j}}
\newcommand{\kkk}{\boldsymbol{k}}
\newcommand{\rrr}{\boldsymbol{r}}
\newcommand{\FFF}{\boldsymbol{F}}
\newcommand{\yyy}{\boldsymbol{y}}
\newcommand{\ppp}{\boldsymbol{p}}
\newcommand{\qqq}{\boldsymbol{q}}
\newcommand{\nnn}{\boldsymbol{n}}
\newcommand{\vvv}{\boldsymbol{v}}
\newcommand{\eee}{\boldsymbol{e}}
\newcommand{\fff}{\boldsymbol{f}}
\newcommand{\www}{\boldsymbol{w}}
\newcommand{\0}{\boldsymbol{0}}
\newcommand{\lon}{\longrightarrow}
\newcommand{\ga}{\gamma}
\newcommand{\pa}{\partial}
\newcommand{\QED}{\hfill $\Box$}
\newcommand{\id}{{\mbox {\rm id}}}
\newcommand{\Ker}{{\mbox {\rm Ker}}}
\newcommand{\grad}{{\mbox {\rm grad}}}
\newcommand{\ind}{{\mbox {\rm ind}}}
\newcommand{\rot}{{\mbox {\rm rot}}}
\newcommand{\diver}{{\mbox {\rm div}}}
\newcommand{\Gr}{{\mbox {\rm Gr}}}
\newcommand{\LG}{{\mbox {\rm LG}}}
\newcommand{\Diff}{{\mbox {\rm Diff}}}
\newcommand{\Symp}{{\mbox {\rm Symp}}}
\newcommand{\Ct}{{\mbox {\rm Ct}}}
\newcommand{\Uns}{{\mbox {\rm Uns}}}
\newcommand{\rank}{{\mbox {\rm rank}}}
\newcommand{\sign}{{\mbox {\rm sign}}}
\newcommand{\Spin}{{\mbox {\rm Spin}}}
\newcommand{\Sp}{{\mbox {\rm sp}}}
\newcommand{\Int}{{\mbox {\rm Int}}}
\newcommand{\Hom}{{\mbox {\rm Hom}}}
\newcommand{\codim}{{\mbox {\rm codim}}}
\newcommand{\ord}{{\mbox {\rm ord}}}
\newcommand{\Iso}{{\mbox {\rm Iso}}}
\newcommand{\corank}{{\mbox {\rm corank}}}
\def\mod{{\mbox {\rm mod}}}
\newcommand{\pt}{{\mbox {\rm pt}}}
\newcommand{\enP}{\hfill $\Box$ \par\vspace{5truemm}}
\newcommand{\qed}{\hfill $\Box$ \par}
\newcommand{\spe}{\vspace{0.4truecm}}
\newcommand{\ad}{{\mbox{\rm ad}}}
%
\newenvironment{FRAME}{\begin{trivlist}\item[]
	\hrule
	\hbox to \linewidth\bgroup
		\advance\linewidth by -10pt
		\hsize=\linewidth
		\vrule\hfill
		\vbox\bgroup
			\vskip5pt
			\def\thempfootnote{\arabic{mpfootnote}}
			\begin{minipage}{\linewidth}}{%
			\end{minipage}\vskip5pt
		\egroup\hfill\vrule
	\egroup\hrule
	\end{trivlist}}

\title{
Duality of singular paths for $(2,3,5)$-distributions
} 

\author{G. Ishikawa\thanks{This work was supported by KAKENHI No.22340030 and No.23654058.}, Y. Kitagawa and W. Yukuno}


\date{ }

\maketitle

\begin{abstract} 
We show a duality which arises from distributions of Cartan type, having growth (2, 3, 5), from the view point of geometric control theory. In fact we consider the space of singular (or abnormal) paths on a given five dimensional space endowed with a Cartan distribution, which form another five dimensional space with a cone structure. We regard the cone structure as a control system and show that the space of singular paths of the cone structure is naturally identified with the original space. Moreover we observe an asymmetry on this duality in terms of singular paths. 
\end{abstract}

\section{Introduction.}

In this paper we show a duality which is related to rank two distributions of five variables: 
Let $D \subset TY$ be a subbundle of the tangent bundle of a five dimensional manifold $Y$
with growth $(2, 3, 5)$ (see \S \ref{Control systems and singular controls.}). 
It is known that for any point $y$ of $Y$ and for any direction of $D_y$,  there exists uniquely a 
singular $D$-path (or an immersed abnormal extremal) 
through $y$ with the given direction (see Proposition \ref{(2,3,5)}). 
Then the immersive 
singular $D$-paths form another five dimensional manifold $X$. 
In fact, the construction of the space $X$ was mentioned in the lecture note \cite{Bryant} shortly. 
Moreover we remark that 
$X$ is endowed with a natural cone structure induced from the prolongation $E$ of $D$ on a six dimensional manifold $Z$
and a projection $Z \to X$. 
In this paper, 
we regard the cone structure $C \subset TX$ as a control system ${\mathbb C} : E \to TX \to X$ on $X$ 
defined by the differential map $E \hookrightarrow TZ \to TX$ of the projection $Z \to X$ and consider the space 
of singular ${\mathbb C}$-paths. 
Then we show that the original space $Y$ is identified with the space of singular ${\mathbb C}$-paths on $X$, 
at least locally, while $X$ is identified with the space of singular $D$-paths (Theorem \ref{duality}). 
In fact, a singular ${\mathbb C}$-path on $X$ 
consists of singular $D$-paths on $Y$ which pass through a fixed point of $Y$. 
%

We recall the basic definitions and terminologies needed in this paper 
in \S \ref{Control systems and singular controls.}, and 
introduce a fundamental lemma for treating singular paths of cone structures 
in \S \ref{Lifting of singular trajectories.}. 
In \S \ref{Duality of singular paths.}, we show our basic constructions and the main duality result 
mentioned above. 

The abnormal extremals of a distribution $D$ with growth $(2, 3, 5)$ 
with $G_2$-symmetry was calculated in \cite{Kitagawa}. 
On the other hand, from the view point of twistor theory, 
the double fibration $Y \leftarrow Z \to X$ with $G_2$-symmetry
and the canonical geometric structures $E \subset TZ, D \subset TY$ and $C \subset TX$ have been explicitly constructed 
in \cite{IMT}. 
Then we can determine singular paths, using the explicit representation, by direct calculations in $G_2$-case. 
In \S \ref{Example: $G_2$-case.}, we recall the explicit representations of the control systems
${\mathbb E} : E \hookrightarrow TZ \to Z$, ${\mathbb D} : D \hookrightarrow TY \to Y$ and 
${\mathbb C} : E \to TX \to X$, and show partly the explicit calculations on their singular paths. 

Distributions with growth $(2, 3, 5)$ was studied by Cartan in \cite{Cartan}. 
In fact they are called systems {\it of Cartan type} in \cite{BH}. 
They are studied in detail from various viewpoints by many authors 
(see for instance \cite{BH}\cite{Zelenko2}\cite{Zhitomirskii}), 
and even now we believe still there remain things to be uncovered on this subject. 

\

All manifolds and mappings are assumed to be of class $C^\infty$ unless otherwise stated. 
For an interval $I$, we say that an assertion holds for almost every $t \in I$ if 
it holds outside of some measure zero set of $I$. 

\section{Control systems and singular controls.}
\label{Control systems and singular controls.}

We need the following generalised notion of control systems (\cite{Agrachev}). 
Let $M$ be a manifold. 
A {\it control system} on $M$ is given by 
a locally trivial fibration $\pi_{\mathcal U} : {\mathcal U} \to M$ over $M$ and 
a fibered map $F : {\mathcal U} \to TM$ over the identity of $M$. 
In this paper, we assume that the fibre $U$ of $\pi_{\mathcal U} : {\mathcal U} \to M$ 
is an open subset of $\R^r$. 
Note that in this paper we mainly treat a local case where $\pi_{\mathcal U}$ is trivial. However 
our results do not depend on any choice of trivialisation of $\pi_{\mathcal U}$. 

For a given control system, an $L^\infty$ (measurable, essentially bounded) 
map $c : [a, b] \to {\mathcal U}$ on an interval is called an {\it admissible control} 
if the curve $\gamma := \pi_{\mathcal U}\circ c : [a, b] \to M$ satisfies 
$\dot{\gamma}(t) = F(c(t))$, for almost every $t \in [a, b]$. 
Then the Lipschitz curve $\gamma$ is called a {\it trajectory}. 
In this paper, we use the term \lq\lq path\rq\rq \ for a smooth ($C^\infty$)
 immersive trajectory regarded up to parametrisation. 

The {\it initial} (resp. {\it end}) {\it point} of $\gamma$ is given by 
$\gamma(a) = \pi_{\mathcal U}(c(a))$ (resp. $\gamma(b) = \pi_{\mathcal U}(c(b))$. 
Note that, in our terminology, the term \lq\lq control" contains 
the data on both the \lq\lq control" in usual sense 
and (the initial point of) the trajectory. Therefore trajectory is determined uniquely by the control, 
while the control is not unique for the trajectory in general. 

If we write locally $x(t) = \gamma(t)$ and $\widetilde{\gamma}(t) = (x(t), u(t))$, 
under a local triviality ${\mathcal U}\vert_V \cong V\times U$ 
over an open set $V \subset M$, 
then the control system is expressed by a family of vector fields, $\{ f_u\}_{u \in U}$, 
$f_u(x) = F(x, u)$ over $V \subset M$ 
and by the equation $\dot{x}(t) = f_{u(t)}(x(t))$. 
In this local situation or in the case where the fibration $\pi_{\mathcal U}$ is globally trivial, 
for a given initial point, $u(t)$ is called a control defining the trajectory $x(t)$ as usual. 
Note that, for a given initial point $x_0$ 
and a given $L^\infty$ control $u(t)$, the Cauchy problem $f_u(x) = F(x, u), x(0) = q_0$, 
has a unique Lipschitz solution (trajectory) $x(t)$ depending smoothly on $q_0$ 
by the classical Carath\'{e}odory theorem (see for example \cite{AS} 2.4.1). 

Now, given a control system and given a point $q_0 \in M$, 
we denote by ${\mathcal C}_{ad}$ 
the set of all admissible controls $c : [a, b] \to {\mathcal U}$ 
with the initial point $\pi_{\mathcal U}(c(a)) = q_0$. 
Then it is known that ${\mathcal C}_{ad}$ is a Banach manifold (\cite{AS}). 
The {\it endpoint mapping} ${\mathcal E} : {\mathcal C}_{ad} \to M$ is defined by 
${\mathcal E}(c) := \pi_{\mathcal U}(c(b))$. 
The control $c$ with the initial point $\pi_{\mathcal U}(c(a)) = q_0$ 
is called {\it singular} or {\it abnormal}, if it is a singular point of ${\mathcal E}$, namely if 
the differential ${\mathcal E}_* : T_c{\mathcal C}_{ad} \to T_{{\mathcal E}(c)}M$ is not surjective.  
If $c$ is a singular control, then the trajectory $\gamma = \pi_{\mathcal U}\circ c$ is called a 
{\it singular trajectory} or an {\it abnormal extremal} (\cite{AS}\cite{BC}). 

We see that if a control is singular then its restriction to any subinterval is singular. 
In fact, as a part of {\it Pontryagin maximal principle}, 
the general characterisation of singular controls is known. 
To state the characterisation, consider the fibre product of  
$\pi_{\mathcal U} : {\mathcal U} \to M$ and $\pi_{T^*M} : T^*M \to M$: 
$$
{\mathcal U} \times_MT^*M := \{ (w, \varphi) \in {\mathcal U}\times T^*M \mid 
\pi_{\mathcal U}(w) = \pi_{T^*M}(\varphi) \}, 
$$
(extended cotangent bundle), endowed with the natural projections $\Pi_{\mathcal U} : {\mathcal U} \times_MT^*M \to {\mathcal U}$ and 
$\Pi_{T^*M} : {\mathcal U} \times_MT^*M \to T^*M$. 
Regarding the local triviality ${\mathcal U}\vert_V \cong V\times U$ 
over a local coordinate neighbourhood $V \subset M$,  
${\mathcal U} \times_MT^*M$ is identified with $({\mathcal U}\vert_V) \times_V (T^*M\vert_V)$ and with 
$T^*M\vert_V\times U$. 

We define the {\it Hamiltonian} function 
$H : {\mathcal U} \times_MT^*M \to \R$ of the control system 
$F : {\mathcal U} \to TM$ by
$$
H(w, \varphi) := \langle \varphi, F(w)\rangle, \quad (w, \varphi) \in {\mathcal U} \times_MT^*M. 
$$
Note that the right-hand side is well defined, since 
$\pi_{T^*M}(\varphi) = \pi_{\mathcal U}(w) = \pi_{TM}(F(w))$. 
In particular, in the case where $U$ is just one point, for a vector field $\xi$ over $M$, we define $H_\xi : T^*M \to \R$ by 
$H_\xi(x, p) = \langle p, \xi(x)\rangle$. 

For each point $\varphi \in T^*M$, 
we set $U_{\varphi} := \Pi_{T^*M}^{-1}(\varphi)$, the fibre of $\Pi_{T^*M} : {\mathcal U} \times_MT^*M \to T^*M$. 
Then $U_{\varphi}$ is diffeomorphic to the fibre $U$ of $\pi_{\mathcal U} : {\mathcal U} \to M$. 

We define the {\it relative critical locus} $\Sigma(H) \subset {\mathcal U} \times_MT^*M$ 
of $H$ by the set of critical points for the restrictions $H\vert_{U_{\varphi}} : 
U_{\varphi} \to \R$, $\varphi \in T^*M$. 

Consider the canonical symplectic form $\omega$ on $T^*M$ 
and the pull-back $\Omega := (\Pi_{T^*M})^*(\omega)$ on ${\mathcal U} \times_MT^*M$ by the mapping $\Pi_{T^*M}$. 
The kernel of the interior product 
$T_{(w, \varphi)}({\mathcal U} \times_MT^*M) \to T^*_{(w, \varphi)}({\mathcal U} \times_MT^*M)$, 
$\xi \mapsto i_{\xi}\Omega$ 
is equal to $T_{(w, \varphi)}U_{\varphi}$, for each point $(w, \varphi) \in {\mathcal U} \times_MT^*M$. 
Then consider the relative exterior differential 
$dH\vert_{\Sigma(H)} : \Sigma(H) \to \Pi_{T^*M}^*(T^*(T^*M))$, 
where $\Pi_{T^*M}^*(T^*(T^*M))$ is the pull-back of the vector bundle $T^*(T^*M)$ over $T^*M$ by 
$\Pi_{T^*M}$. 
It defines a vector field 
$\overrightarrow{H} : \Sigma(H) \to \Pi_{T^*M}^*(T(T^*M))$ along $\Sigma(H)$, by 
$$
i_{\overrightarrow{H}}{\Omega}(w, \varphi) = -dH(w, \varphi), \quad (w, \varphi) \in \Sigma(H), 
$$
We call $\overrightarrow{H}$ the {\it Hamilton vector field}, 
which is just the family of Hamilton vector fields with parameter $u \in U$ restricted to $\Sigma(H)$, 
in local situation over $M$. 
Let $O \subset T^*M$ denote the zero-section of $\pi_{T^*M} : T^*M \to M$. 
Then we have: 

\bep {\rm (\cite{AS})} 
\label{PMP}
A control $c : [a, b] \to {\mathcal U}$ is a singular control if and only if 
there exists a lift $\beta : [a, b] \to ({\mathcal U}\times_MT^*M) \setminus \Pi_{T^*M}^{-1}(O)$ of $c$ for 
$\Pi_{\mathcal U} : {\mathcal U}\times_MT^*M \to {\mathcal U}$ which 
satisfies that $\Pi_{T^*M}\circ \beta : [a, b] \to T^*M$ is Lipschitz and satisfies 
the constrained Hamiltonian equation: 
$$
(\Pi_{T^*M}\circ \beta)'(t) = \overrightarrow{H}(\beta(t)), \quad \beta(t) \in \Sigma(H) (\subset {\mathcal U}\times_MT^*M), 
$$
for almost every $t \in [a, b]$. 
\enp

The lift in the above Proposition \ref{PMP} is called {\it a singular bi-extremal} or {\it an abnormal bi-extremal} for $c$. 

In particular the notion of singular controls is a local notion: 

\bec
A control $c : I \to {\mathcal U}$ is a singular control 
if and only if, for any subinterval $J \subset I$, the restriction 
$c : J \to {\mathcal U}$ of $c$ is a singular control. 
\enc

The local description of the characterisation of singular control is given as follows: 

\bep
Suppose $\dim M = m, \dim U = r$. 
Let 
$$
(x; p; u) = (x_1, \dots, x_m; p_1, \dots, p_m; u_1, \dots, u_r)
$$ 
be local coordinates of ${\mathcal U}\vert_V\times_VT^*M\vert_V \cong T^*M\vert_V \times U$ 
over coordinate neighbourhood $V \subset M$ and 
$H(x; p; u)$ the local expression of the Hamiltonian of a given control system. 
Then a control $c : [a, b] \to V\times U$ is singular if and only if 
there exists a lift $\beta : [a, b] \to T^*M\vert_V \times U$ of $c$, 
$\beta(t) = (x(t); p(t); u(t))$, which satisfies, for almost every $t \in [a, b]$, 
$$
\begin{cases}
\ \dot{x}_i(t)  =  \ \ \dfrac{\pa H}{\pa p_i}(x(t); p(t); u(t)), \quad (1 \leq i \leq m)
\vspace{0.2truecm}
\\
\ \dot{p}_i(t)  =  - \dfrac{\pa H}{\pa x_i}(x(t); p(t); u(t)), \quad (1 \leq i \leq m)
\vspace{0.2truecm}
\\
\ \dfrac{\pa H}{\pa u_j}(x(t); p(t); u(t)) = 0, \quad (1 \leq j \leq r), \qquad 
p(t) \not= 0. 
\end{cases}
$$
\enp

Consider two control systems $F, G : {\mathcal U} \to TM$ 
on $M$. Given a local trivialisation ${\mathcal U}\vert_V \cong V \times U$ over an open subset 
$V \subset M$, 
regard $F, G$ as families of vector fields $f_u, g_u$ respectively. Then we define a control system 
$[F, G] : {\mathcal U}\vert_V \to TV$ by $[F, G](x, u) = [f_u, g_u](x), x \in V, u \in U$. 
We denote by $H_F$ the Hamiltonian function $H : T^*V\times U \to \R$ when we make stress on a control system $F$. 
In the next section, we will use the following fundamental formula (cf. \cite{AS}, Chapter 11, for instance): 

\bel
\label{Hamiltonian-general}
We have 
$dH_G(\overrightarrow{H}_F) = H_{[F, G]} : T^*V \times U \to \R$. 
\enl

In particular we have, in the case where $U = \pt$, 

\bel
\label{Hamiltonian}
For vector fields $\xi, \eta$ over $M$, 
we have $dH_\xi(\overrightarrow{H}_\eta) = H_{[\eta, \xi]} : T^*M \to \R$. 
\enl

\

A vector subbundle $D \subset TM$ is called a {\it distribution} or a {\it differential system} on $M$. Then 
the inclusion mapping $F : D \to TM$ with the bundle projection $\pi : D \to M$ defines a 
control system in the above sense. 
In this case, a control (resp. a trajectory, a path) 
is called a $D$-control (resp. a $D$-trajectory, a $D$-path). 
Then a $D$-control is uniquely determined from its $D$-trajectory. 
Let $\xi_1, \dots, \xi_r$ be a local frame of $D$. Then the Hamiltonian is given by 
$$
H_{\xi_i}(x; p) := \langle p, \xi_i(x)\rangle
\ \ {\mbox{\rm and }} \ 
H(x; p; u) := u_1H_{\xi_1}(x; p) + \cdots + u_rH_{\xi_r}(x; p). 
$$
The local characterisation of singular $D$-controls is provided by the existence of $p(t) \not= 0$ 
satisfying 
$$
\begin{cases}
\ 
\dot{x}(t) = \ \ u_1\xi_1(x(t)) + \cdots + u_r\xi_r(x(t)), 
\vspace{0.2truecm}
\\
\ 
\dot{p_i}(t) = 
- u_1\dfrac{\ \pa H_{\xi_1}}{\pa x_i}(x(t); p(t)) - \cdots - u_r\dfrac{\
 \pa H_{\xi_r}}{\pa x_i}(x(t); p(t)), 
\quad(1 \leq i \leq m), 
\vspace{0.2truecm}
\\
\ 
H_{\xi_i}(x(t); p(t)) = 0, \quad (1 \leq j \leq m). 
\end{cases}
$$

\

Let $D \subset TM$ be a distribution and ${\mathcal D}$ the sheaf of section-germs to $D$. 
The {\it small} (or {\it weak}) (resp. the {\it big} (or {\it strong})) 
derived systems ${\mathcal D}^{(i)}$, (resp. ${\mathcal D}^i$), $i = 1, 2, 3, \dots, $ are defined by 
${\mathcal D}^{(1)} = {\mathcal D}^1 = {\mathcal D}$ and, inductively, by 
$$
{\mathcal D}^{(i+1)} := {\mathcal D}^{(i)} + [{\mathcal D}, {\mathcal D}^{(i)}], 
 \ \ ({\mbox{\rm resp. }} \
{\mathcal D}^{i+1} := {\mathcal D}^{i} + [{\mathcal D}^{i}, {\mathcal D}^{i}]), \ i \geq 1.  
$$
Note that ${\mathcal D}^{(2)} = {\mathcal D}^{2}$ and ${\mathcal D}^{(i)} \subset {\mathcal D}^{i}, 
\ i \geq 3$. 
The distribution $D$ is called {\it bracket generating} if 
$\cup_{i = 1}^\infty {\mathcal D}^{(i)} = {\mathcal T}M$, the total sheaf of vector fields over $M$. 
Suppose each ${\mathcal D}^{(i)}$ (resp. ${\mathcal D}^{i}$) coincides with the sheaf of section-germs to some subbundle of 
$TM$, written as $D^{(i)}$ (resp. $D^{i}$). In this situation, we call $D$ has {\it small} (or {\it weak}) (resp. {\it big} (or {\it strong})) 
growth $(n_1, n_2, \dots, n_i, \dots)$ if $\rank(D^{(i)}) = n_i$ (resp. $\rank(D^{i}) = n_i), \ i \geq 1$. 

A singular control $c : I \to D$ of the control system $D \hookrightarrow TM \to M$ is called 
{\it regular} 
if there exists a singular bi-extremal $\beta : I \to D \times_M T^*M$ for $c$ such that 
$\Pi_{T^*M}\circ \beta(t) \in D^{(2)\perp} \setminus D^{(3)\perp} \subset T^*M$ 
for any $t \in I$ (\cite{Montgomery}\cite{Hsu}\cite{LS}). 
A singular control $c : I \to D$ is called 
{\it totally irregular} if any singular bi-extremal $\beta : I \to D \times_M T^*M$ for $c$ satisfies 
that $\Pi_{T^*M}\circ \beta(t) \in D^{(3)\perp} \subset T^*M$ 
for any $t \in I$. Here, for a subbundle $E \subset TM$, we set 
$E^{\perp} = \{ p \in T^*M \mid \langle p, v \rangle = 0, {\mbox{\rm \ for any \ }} 
v \in E_{\pi_{T^*M}(p)} \}$. 

A singular trajectory (resp. a singular path) is called {\it regular} (resp. {\it totally irregular}) 
if so is its control. Note again that in the case of distribution, 
a control $c$ is uniquely determined from the trajectory $\pi_D\circ c$.

\section{Liftings of singular trajectories.}
\label{Lifting of singular trajectories.}

Let $M, N$ be manifolds of dimension $m = n + k, n$ respectively, 
$E \subset TM$ a distribution on $M$ of rank $r = \ell + k$, 
and $\pi : M \to N$ a fibration with $K = \Ker(\pi_*) \subset E$. We only treat the case where 
the $\pi$-fibre is an open subset of $\R^k$ and $E$ is trivial over on each $\pi$-fibre. 

Consider two control systems 
$$
{\mathbb E} : E \hookrightarrow TM \xrightarrow{\pi_{TM}} M, 
{\mbox{\rm \  over \ }} M, \quad {\mbox{\rm \  and \ }}\quad
{\mathbb E}/\pi : E \xrightarrow{\pi_*\vert_E} TN \xrightarrow{\pi_{TN}} N, 
{\mbox{\rm \ over \ }} N, 
$$
respectively. 
Note that the composition $\pi\circ \pi_E : E \to N$ of the bundle projection $\pi_E : E \to M$ 
and the fibration $\pi : M \to N$ is again a fibration over $N$, having the fibre $\pi_E^{-1}(W)$ for a 
typical fibre $W \subset M$ of $\pi$. 

Moreover we suppose $E = L \oplus K$ for a complementary subbundle $L \subset E$ to $K$ in $E$ of rank $\ell$ and 
consider the third control system 
$$
{\mathbb L} : L \xrightarrow{\pi_*\vert_L} TN \xrightarrow{\pi_{TN}} N, 
{\mbox{\rm \ over \ }} N. 
$$

Then first we have 

\bel
\label{lifting}
Let $\gamma : I \to N$ be a singular $({\mathbb E}/\pi)$-trajectory. Suppose that 
there exists a Lipschitz abnormal bi-extremal $\beta : I \to E\times_N T^*N$ corresponding to $\gamma$. 
Then there exists  a singular ${\mathbb E}$-trajectory $\widetilde{\gamma} : 
I \to M$ such that $\pi\circ \widetilde{\gamma} = \gamma$. 
\enl

\ber
{\rm
In Lemma \ref{lifting}, 
we pose the Lipschitz condition on abnormal bi-extremals, because we have to regard some of control parameters 
as state variables. 
}
\enr

\noindent
{\it Proof of Lemma \ref{lifting}: }
It is sufficient to deal with local situation on $N$: 
Take a system of local coordinates 
$x_1, \dots, x_n, w_1, \dots, w_k$ on $M$ such that 
$\pi(x, w) = x$ and therefore $K = \langle \frac{\pa}{\pa w_1}, \dots, \frac{\pa}{\pa w_k} \rangle$. 
Take a local frame $\xi_1, \dots, \xi_\ell$ of $L$, 
$
\xi_i(x, w) = \sum_{j=1}^n c_{ij}(x, w)\frac{\pa}{\pa x_j}, (1 \leq i \leq \ell), $ so that 
$\xi_1, \dots, \xi_\ell, \frac{\pa}{\pa w_1}, \dots, \frac{\pa}{\pa w_k}$ form a local frame of $E$. 
The system of canonical local coordinates of $T^*M$ is given by 
$(x, w; p_1, \dots, p_n, \psi_1, \dots, \psi_k)$. 
For instance, $K^{\perp} \subset T^*M$ is defined locally by $\psi_1 = 0, \dots, \psi_k = 0$. 
For the control system ${\mathbb E}$, the Hamiltonian $H$ is given by 
$$
H(x, w; p, \psi; \lambda, \mu) = \sum_{1 \leq i \leq \ell} \lambda_i H_{\xi_i}(x, w; p, \psi) 
+ \sum_{1 \leq i \leq k} \mu_iH_{\pa/\pa w_i}(x, w; p, \psi), 
$$
where 
$H_{\xi_i}(x, w; p, \psi) = \sum_{1\leq j \leq n} c_{ij}(x, w)p_j$, $H_{\pa/\pa w_i}(x, w; p, \psi) = \psi_i$, and 
$\lambda = (\lambda_1, \dots, \lambda_\ell)$ (resp, $\mu = (\mu_1, \dots, \mu_k)$) are the fibre coordinates of $L$ (resp. $K$). 
The constrained Hamiltonian system is given by 
$$
\begin{array}{c}
{\displaystyle 
\dot{x}_j = \sum_{1\leq i \leq \ell} \lambda_i c_{ij}, (1 \leq j \leq n), \quad \dot{w}_j = \mu_j, (1 \leq j \leq k), }
\vspace{0.2truecm}
\\
{\displaystyle 
\dot{p}_{\nu} = - \sum_{1\leq i \leq \ell, 1 \leq j \leq n}\lambda_i\dfrac{\pa c_{ij}}{\pa x_{\nu}} p_j, 
(1 \leq \nu \leq n), \quad 
\dot{\psi}_{\kappa} = - \sum_{1\leq i \leq \ell, 1 \leq j \leq n} \lambda_i\dfrac{\pa c_{ij}}{\pa w_{\kappa}} p_j, 
(1 \leq \kappa \leq k), } 
\vspace{0.2truecm}
\\
{\displaystyle 
\sum_{1\leq j \leq n} c_{ij}p_j = 0, (1 \leq i \leq \ell), \quad \psi_i = 0, (1 \leq i \leq k). 
}
\end{array}
$$

The control system ${\mathbb E}/\pi$ over $N$ is locally described as follows: 
The system of local coordinates of $N$ (resp. $T^*N$, $E$) 
is given by $x = (x_1, \dots, x_n)$, (resp. $(x; p) = (x; p_1, \dots, p_n)$, $(x, w; \lambda, \mu) 
= (x, w;\lambda_1, \dots, \lambda_{\ell}, \mu_1, \dots, \mu_k)$). Then the system of local coordinates of 
$E\times_N T^*N$ is given by $(x, w; p; \lambda, \mu)$. Note that $\dim(E\times_N T^*N) = 2n + 2k + \ell$. 
In this case, the mapping $F : E \to TN$ is given locally by 
$$
F(x, w; \lambda, \mu) = \sum_{1 \leq i \leq \ell} \lambda_i\, \pi_{X*}\xi_i(x, w) + 
\sum_{1 \leq i \leq n} \mu_i\, \pi_{X*}\frac{\pa}{\pa w_i}(x, w) 
= \sum_{1 \leq i \leq \ell} \lambda_i\, \pi_{X*}\xi_i(x, w), 
$$
and the Hamiltonian function of the control system ${\mathbb E}/\pi$ is given by 
$$
H(x, w; p; \lambda, \mu) = \langle p, \ \sum_{1 \leq i \leq \ell} \lambda_i\, \pi_{X*}\xi_i(x, w) \rangle 
= \sum_{1 \leq i \leq \ell} \lambda_i \langle p, \pi_{X*}\xi_i(x, w)\rangle
= \sum_{1\leq i \leq \ell, 1 \leq j \leq n} \lambda_i c_{ij} p_j. 
$$ 
Here the control parameters 
are given by $w, \lambda, \mu$. 
Note that $H$ is independent of $\mu$. The constrained Hamiltonian system is given 
by 
$$
\begin{array}{c}
{\displaystyle 
\dot{x}_j = \sum_{1\leq i \leq \ell} \lambda_i c_{ij}, (1 \leq j \leq n), \quad 
\dot{p}_{\nu} = - \sum_{1\leq i \leq \ell, 1 \leq j \leq n} \lambda_i\dfrac{\pa c_{ij}}{\pa x_{\nu}} p_j, 
(1 \leq \nu \leq n), 
}
\vspace{0.2truecm}
\\
{\displaystyle 
\sum_{1\leq j \leq n} c_{ij}p_j = 0, (1 \leq i \leq \ell), \quad 
\sum_{1\leq i \leq \ell, 1 \leq j \leq n} \lambda_i\dfrac{\pa c_{ij}}{\pa w_{\kappa}} p_j = 0, 
(1 \leq \kappa \leq k). 
}
\end{array}
$$

Let $\beta(t) = (x(t), w(t); p(t); \lambda(t), \mu(t))$ be an abnormal bi-extremal for ${\mathbb E}/\pi$. 
Note that $\mu(t)$ can be taken as arbitrary $L^\infty$ function. 
Suppose $\beta(t)$ is Lipschitz. Then in particular the control $w(t)$ is Lipschitz. 
Replace $\mu(t)$ by $\dot{w}(t)$, which is of class $L^\infty$, and take $\psi(t) = 0$.
We set $\widetilde{\beta}(t) = (x(t), w(t); p(t), 0; \lambda(t), \dot{w}(t))$. Then $\widetilde{\beta}$ is 
an abnormal bi-extremal for ${\mathbb E}$ and $\widetilde{\gamma}(t) = (x(t), w(t))$ 
is a lift of $\gamma(t) = x(t), \pi\circ \widetilde{\gamma} = \gamma$. 
\QED

\

Second we remark that 
the projection $\pi_L : E = L \oplus K \to L$ induces the  
projection $\rho : E\times_N T^*N \to L \times_N T^*N$. Then we have 

\bel
\label{lifting2}
A curve $\beta : I \to E \times_N T^*N$ is an abnormal bi-extremal (resp. 
a Lipschitz abnormal bi-extremal) for ${\mathbb E}/\pi$ if and only if 
$\rho\circ\beta : I \to L \times_N T^*N$ is an abnormal bi-extremal (resp. 
a Lipschitz abnormal bi-extremal) for ${\mathbb L}$. Moreover any abnormal bi-extremal
(resp. a Lipschitz abnormal bi-extremal) $\overline{\beta} : I \to L \times_N T^*N$ for ${\mathbb L}$ 
is written as $\rho\circ\beta$ by an abnormal bi-extremal (resp. 
a Lipschitz abnormal bi-extremal) $\beta : I \to E \times_N T^*N$ for ${\mathbb E}/\pi$. 
\enl

\Proof
The system of local coordinates of $L \times_N T^*N$ is given 
$(x, w;p ; \lambda)$, using the coordinates in the proof of Lemma \ref{lifting}. 
The constrained Hamiltonian system for ${\mathbb L}$ 
is of the same form for ${\mathbb E}/\pi$ if the un-efficient component $\mu$ is deleted. 
Therefore $\beta(t) = (x(t), w(t); p(t); \lambda(t), \mu(t))$ is an abnormal bi-extremal (resp. 
a Lipschitz abnormal bi-extremal) for ${\mathbb E}/\pi$ if and only if 
$\rho\circ\beta(t) = (x(t), w(t); p(t); \lambda(t))$ is an abnormal bi-extremal (resp. 
a Lipschitz abnormal bi-extremal) for $L$. 
The second assertion is clear. 
\QED

\

Let $z \in M$. Let $\xi$ be a section-germ at $z$ to $L$ and $\kappa$ a section-germ at $z$ to $K$. 
Then $[\xi, \kappa](z)$ modulo $E_z$ depends only on the tangent vectors $\xi(z), \kappa(z) \in T_zM$. 
For $v \in L_z$, we define a linear map $\ad(v) : K_z \to T_zM/E_z$ by 
$\ad(v)(u) := [\xi, \kappa](z) \  \mod\, E_z, \xi(z) = v, \kappa(z) = u$. 
We need the following result for the proof of our main result in the next section. 

\bel
\label{immersion-general}
Let $(z; v) \in L, v \in L_z$. Then 
the differential map 
$\pi_{X*} \vert_L : L  \to TX$ is an immersion at $(z; v)$ 
if and only if $\ad(v) : K_z \to T_zM/E_z$ is injective. 
\enl

\Proof
We use the system of coordinates $(x, w)$ of $M$ and $(x, w; \lambda)$ of $L$ 
as in Lemma \ref{lifting} and, take a frame of $L$, 
$\xi_i(x, w) = \sum\limits_{1\leq j \leq n} c_{ij}(x, w)\frac{\pa}{\pa x_j}, (1 \leq i \leq \ell). $
Then $\pi_* : L \to TX$ is expressed by 
$$
(x, w; \lambda) \mapsto (x_1, \dots, x_n, \sum_{1\leq i \leq \ell} \lambda_ic_{i1}, \dots, \sum_{1\leq i \leq \ell} \lambda_ic_{in}). 
$$
The Jacobi matrix of the mapping $\pi_*$ is given by 
$$
\left(
\begin{array}{ccccccc}
E_n & O & \cdots & O & O & \cdots & O \\
* & \sum\limits_{i} \lambda_i\frac{\pa c_{i1}}{\pa w_1} & \cdots & \sum\limits_{i} \lambda_i\frac{\pa c_{i1}}{\pa w_k} & c_{11} & \cdots & 
c_{\ell 1} 
\\
* & \vdots & \ddots & \vdots & \vdots & \ddots & \vdots 
\\
* & \sum\limits_{i} \lambda_i\frac{\pa c_{in}}{\pa w_1} & \cdots & \sum\limits_{i} \lambda_i\frac{\pa c_{in}}{\pa w_k} & c_{1n} & \cdots & 
c_{\ell n} 
\end{array}
\right). 
$$
Let $(z; v) \in L$ and $v = \sum\limits_{1\leq i \leq \ell} \lambda_i\xi_i(z)$. 
Then the Jacobi matrix at $(z; v) \in L$ is of rank $n + k + \ell$ if and only if 
$$
{\textstyle 
\left[\sum\limits_{1\leq i \leq \ell} \lambda_i\xi_i, \ \frac{\pa}{\pa w_1}\right], \dots, 
\left[\sum\limits_{1\leq i \leq \ell} \lambda_i\xi_i, \ \frac{\pa}{\pa w_k}\right], \xi_1, \dots, \xi_{\ell}
}
$$
are linear independent at $z$, which is equivalent to that $\ad(v)$ is injective. 
\QED


\section{Duality of singular paths.}
\label{Duality of singular paths.}

Let $Y$ be a $5$-dimensional manifold and $D \subset TY$ a distribution on $Y$ of 
rank two. 
Let $Z$ be the set of tangential directions on $D$, 
$$
Z : = \{ z \mid z  {\mbox{\rm \ is an oriented one-dimensional linear subspace of 
}} 
D_y 
{\mbox{\rm \ for some\ }} y \in Y \}, 
$$
with the natural projection $\pi_Y : Z \to Y$. 
Define the prolongation $E \subset TZ$ of $D$, tautologically from $D$ by
$E_z := (\pi_{Y*})^{-1}(z)$ for any $z \in Z$ regarded as a line $z \subset T_{\pi_Y(z)}Y$ (\cite{BH}). 
Suppose $D$ has small growth $(2, 3, 5)$. 
Note that $D$ has big growth $(2, 3, 5)$ as well. 
We give a detailed proof of the following known result, in order to use in the process of the proof on 
our main result. 

\bep 
\label{(2,3,5)}
{\rm (\cite{LS}\cite{Montgomery})}
Let $D \subset TY$ be a distribution with growth $(2, 3, 5)$ on a $5$-dimensional manifold 
$Y$. Then, given any point $y \in Y$ and any oriented $\ell \subset T_yY$, 
there exists an oriented singular $D$-trajectory through $p$ with 
the direction $\ell$, which is smooth and is unique up to parametrisation. 
\enp

\Proof
Let $D$ be a distribution with growth $(2, 3, 5)$ 
on a $5$-dimensional manifold $Y$. Recall that $D$ is of rank $2$, $D^2$ is of rank 
$3$ and $D^{(3)} = TY$. 
Let $\{\eta_1, \eta_2\}$ be a local frame of $D$. 
Then $\eta_1, \eta_2$ and $[\eta_1, \eta_2] = \eta_3$ form a 
local frame of $D^2$. Therefore, locally, say on an open set $U \subset Y$, there exists 
a frame $\{\eta_1, \eta_2, \eta_3, \eta_4, \eta_5\}$ of $TU$ such that 
$$
D = \langle \eta_1, \eta_2\rangle, 
\ \ [\eta_1, \eta_2] = \eta_3, \ \ [\eta_1, \eta_3] = \eta_4, \ \ 
[\eta_2, \eta_3] = \eta_5. 
$$

Then the Hamiltonian function $H$ on $D \times_Y T^*Y = T^*U\times \R^2$ is given by 
$$
H(y, q; u) = u_1H_{\eta_1}(y; q) + u_2H_{\eta_2}(y; q), 
$$
where $H_{\eta_i}(y; q) = \langle q, \eta_i(y)\rangle, \ i = 1, 2$, 
and the constrained Hamiltonian system is given by 
$$
\begin{cases}
\ 
\dot{y}(t) = u_1(t)\eta_1(y(t)) + u_2(t)\eta_2(y(t)), 
\vspace{0.2truecm}
\\
\ 
\dot{q}_i(t) = - u_1(t)\dfrac{\ \pa H_{\eta_1}}{\pa y_i}(y(t); q(t))
- u_2(t)\dfrac{\ \pa H_{\eta_2}}{\pa y_i}(y(t); q(t)), \quad (1 \leq i \leq 5), 
\vspace{0.2truecm}
\\
\ 
H_{\eta_1}(y(t); q(t)) = 0, \quad H_{\eta_2}(y(t); q(t))  = 0, \quad q(t) \not= 0. 
\end{cases}
$$
Differentiating both sides of $H_{\eta_1}(y(t); q(t)) (= \langle q(t), \eta_1(x(t))\rangle) = 0$ by $t$, 
we have, using Lemma \ref{Hamiltonian-general}, 
$$
\begin{array}{rcl}
0 & = & (dH_{\eta_1}(\overrightarrow{H}))(y(t); q(t); u(t)) = 
dH_{\eta_1}(u_1\overrightarrow{H}_{\eta_1} +
u_2\overrightarrow{H}_{\eta_2})
(y(t); q(t); u(t))
\\
 & = & u_2(t)H_{[\eta_2, \eta_1]}(y(t); q(t)). 
 \end{array}
$$
Similarly we have, from $H_{\eta_2}(y(t); q(t)) = 0$, 
that $u_1(t)H_{[\eta_1, \eta_2]}(y(t); q(t)) = 0$. 
Therefore if $(u_1(t), u_2(t)) \not= (0, 0)$, then $H_{\eta_3}(y(t); q(t)) = 0$. 
Thus we have that a Lipschitz curve $y(t)$ is a locally non-constant singular $D$-trajectory 
if and only if there exist $u(t)$ and $q(t) \in (D^2)^\perp \setminus 0$ 
satisfying the above constrained Hamiltonian system with 
$H_{\eta_3}(y(t); q(t)) = 0$ in addition. 
Moreover from $H_{\eta_3}(y(t); q(t)) = 0$, we have 
$$
u_1(t) H_{\eta_4}(y(t); q(t)) + u_2(t) H_{\eta_5}(y(t); q(t)) = 0. 
$$
Consider the vector field 
$$
{Ab}_{\eta_1, \eta_2} := H_{\eta_5}\overrightarrow{H}_{\eta_1} - H_{\eta_4}\overrightarrow{H}_{\eta_2}
$$
over $T^*Y$ (see \cite{Zelenko1}). Then ${Ab}_{\eta_1, \eta_2}$ is tangent to $(D^2)^\perp$. 
In fact 
$$
{Ab}_{\eta_1, \eta_2}(H_{\eta_1}) = 0, {Ab}_{\eta_1, \eta_2}(H_{\eta_2}) = 0 {\mbox{\rm \  \ and \ }}
{Ab}_{\eta_1, \eta_2}(H_{\eta_3}) = 0, 
$$ 
on $(D^2)^\perp \subset T^*Y$. 
Therefore 
any singular (oriented) path is obtained as the projection of an integral curve of ${Ab}_{\eta_1, \eta_2}$ 
or $-{Ab}_{\eta_1, \eta_2}$ in $(D^2)^\perp$. 
\qed

\ 

Let $[u_1, u_2]$ be a projective coordinate of the fibre of $\pi_Y : Z = PD \to Y$ and 
$z = u_2/u_1$ an affine coordinate of the fibre. 
Locally identify $Z$ with $\R\times Y$. 
Set $\zeta = \pa/\pa z$. 
Then the prolongation $E$ is generated by $\zeta$ and $\eta_1 + z\eta_2$.

\bel
\label{prolongation}
The prolongation of a distribution with growth $(2, 3, 5)$ has small growth $(2, 3, 4, 5, 6)$ and 
big growth $(2, 3, 4, 6)$. 
\enl

\Proof
Since 
$
[\zeta, \eta_1 + z\eta_2] = [\zeta, \eta_1] + [\zeta, z\eta_2] = \eta_2, 
$
we have 
$
E^2 = \langle \zeta, \eta_1 + z\eta_2, \eta_2\rangle = \langle \zeta, \eta_1, \eta_2\rangle, 
$
which is of rank $3$. 
Since $[\eta_1 + z\eta_2, \eta_2] = \eta_3$, we have 
$E^{(3)} = E^3 = \langle \zeta, \eta_1, \eta_2, \eta_3\rangle$, which is of rank $4$. 
Since $[\eta_1 + z\eta_2, \eta_3] = \eta_4 + z\eta_5$, we have 
$E^{(4)} = \langle \zeta, \eta_1, \eta_2, \eta_3, \eta_4 + z\eta_5\rangle$, that is of rank $5$, while 
$E^4 = TZ$, that is of rank $6$. 
Since $[\zeta, \eta_4 + z\eta_5] = \eta_5$, we have $E^{(5)} = TZ$. 
\qed

\

Note that, by Proposition \ref{(2,3,5)}, any immersed $D$-trajectory lifts to a $E$-trajectory, and 
$Z$ is foliated by lifted paths of singular $D$-paths. 
Thus we have a subbundle $K \subset E$ by this foliation. 
We set $L = \Ker(\pi_{Y*})$. Then $E = L \oplus K$. 
As we see below, this decomposition is intrinsically determined from $(Z, E)$. 

The notions of regular singular paths and totally irregular paths are defined at the end of 
\S \ref{Control systems and singular controls.}. Then we have: 


\bep
\label{E-singular}
Let $D$ be a 
distribution with growth $(2, 3, 5)$ on a five dimensional manifold $Y$, and $(Z, E)$ the prolongation of $D$. 
Then any singular $E$-path (un-parametrised immersed $E$-trajectory) 
is either a fibre of $\pi_Y : Z \to Y$ or the lift of a singular $D$-path. 
Any $\pi_Y$-fibre is a regular singular path of $E$. 
The lift of any singular $D$-path is a totally irregular singular path of $E$. 
\enp

\Proof
Consider the Hamiltonian of $E$: 
$$
H(y, z; q, \varphi; \lambda, \mu) = \lambda\varphi + \mu H_{\eta_1 + z\eta_2}(y, z; q, \varphi), 
$$
where $(y, z; q, \varphi)$ is the canonical coordinate system of $T^*Z$ 
for a coordinate system $(y, z)$ of $Z$ and $\lambda, \mu$ are control parameters. 
Then the constrained Hamiltonian equation for an abnormal bi-extremal 
$$
\beta(t) = (y(t), z(t); q(t), \varphi(t); \lambda(t), \mu(t))
$$ is given by 
$$
\left\{ 
\begin{array}{rcl}
\dot{z}(t) & = & \lambda(t)
\\
\dot{y}(t) & = & \mu(t)(\eta_1 + z\eta_2)(y(t), z(t))
\\
\dot{\varphi}(t) & = & - \mu(t)H_{\eta_2}(y(t); q(t))
\\
\dot{q}(t) & = & - \mu(t)\frac{\pa (H_{\eta_1 + z\eta_2})}{\pa y}(y(t), z(t); q(t), \varphi(t))
\end{array}
\right.
$$
with
$\varphi(t) = 0, H_{\eta_1 + z\eta_2}(y(t), z(t); q(t), \varphi(t)) = 0, (q(t), \varphi(t)) \not= 0$. 

We are supposing $(y(t), z(t))$ is $C^\infty$ and $(\dot{y}(t), \dot{z}(t)) \not= 0$ for all $t \in I$. 
Therefore $(\lambda(t), \mu(t))$ is $C^\infty$ and $(\lambda(t), \mu(t)) \not= 0$. 
Suppose there is a non-empty sub-interval $J \subset I$ where $\mu(t) \not= 0$. Then, on $J$, we have 
$$
H_{\eta_1}(y(t); q(t)) = 0, \quad H_{\eta_2}(y(t); q(t)) = 0. 
$$
Then we have 
$$
0 = \lambda(t)H_{[\zeta, \eta_2]}(y(t); q(t))  + \mu(t)H_{[\eta_1 + z\eta_2, \eta_2]}(y(t); q(t)) 
= \mu(t)H_{\eta_3}(y(t); q(t)). 
$$
Since $\mu(t) \not= 0$, we have $H_{\eta_3}(y(t); q(t)) = 0$. Thus $q(t) \in D^2$ and 
$(q(t), 0) \in E^{(3)\perp}$. Moreover we have 
$$
H_{\eta_4}(y(t); q(t)) + z(t)H_{\eta_5}(y(t); q(t)) = 0. 
$$
Therefore, on $J$, the immersive singular $E$-trajectory is the lift of an immersive singular $D$-trajectory and 
$(q(t), 0) \in E^{(4)\perp}$. 
Take $J$ maximal. Then $J = I$. In fact, otherwise, take an boundary point $t_0$ 
of $J$ which is an interior point of $I$. Then $\lambda (t_0) \not= 0$ and $\mu(t_0) = 0$, so 
$(\dot{y}(t_0), \dot{z}(t_0)) = (0, \lambda(t_0))$ must be tangent to the lift of an immersive singular $D$-trajectory,  
which is a contradiction. 
Thus, in this case, we have that the trajectory $(y(t), z(t))$ coincides with 
the lift of an immersive singular $D$-trajectory on $I$, which is totally irregular. 

If $\mu(t) = 0$ on $I$, then, $y(t)$ is constant on $I$. Setting $y(t) = y_0$,  
we have that the trajectory $(y_0, z(t))$ gives a $\pi_Y$-fibre. 
In this case, for any $q_0 \not= 0$ with 
$H_{\eta_1}(y_0; q_0) = 0, H_{\eta_2}(y_0; q_0) = 0, 
H_{\eta_3}(y_0; q_0) \not= 0$, we have 
a singular bi-extremal $(y_0, z(t); q_0, 0; \dot{z}(t), 0)$ and 
$(q_0, 0) \in E^{2\perp} \setminus E^{(3)\perp}$. Therefore any $\pi_Y$-fibre is 
regular singular. 
\QED

\

\ber
{\rm 
The lift of a singular $D$-path is irregular singular of \lq\lq parabolic type\rq\rq \ in the sense of \cite{Zelenko1}. 
}
\enr

Thus we have locally a fibration $\pi_X : Z \to X$ to another $5$-dimensional manifold $X$ 
and thus the double fibration 
$$
Y \xleftarrow{\ \pi_Y\ } Z \xrightarrow{\ \pi_X\ } X. 
$$
In fact the construction is local in $Y$: If $Y$ is a germ at a point $y_0 \in Y$, 
then we can take $Z$ as a germ along the simple closed curve $\pi_Y^{-1}(y_0) \cong S^1$ and 
$X$ as a germ along the simple closed curve $\pi_X(\pi_Y^{-1}(y_0)) \cong S^1$. 
The prolongation $E \subset TZ$ is expressed as $\Ker(\pi_{Y*}) \oplus \Ker(\pi_{X*})$. 
Therefore, for each $z \in Z$, $E_z$ projects to a line $\ell_z := \pi_{X*}(E_z) \subset T_{\pi_X(z)}X$ 
by $\pi_{X*} : T_zZ \to T_{\pi_X(z)}X$. 
Thus, for any $x \in X$, we have a family of lines 
$\{ \ell_z \mid z \in \pi_X^{-1}(x)\} \subset T_{x}X$. 
We define a cone
$$
C_x := \bigcup _{z \in \pi_X^{-1}(x)} \ell_z \ \subset  T_{x}X, 
$$
and the cone field $C \subset TX$. 
This cone field is regarded as a control system as follows: 
Let $L = \Ker(\pi_{Y*}) \subset E$. Note that $\ell_z = \pi_{X*}(L_z)$. 
We set $F := \pi_{X*}\vert_L : L \to TX$, the composition of the inclusion to $TZ$ and the differential 
$\pi_{X*} : TZ \to TX$ of the projection $\pi_X$. The projection is defined by $\pi_{L} = \pi_{TX}\circ F : L \to X$. 
We denote this control system by ${\mathbb C}$. Note that, in other word,  
${\mathbb C}$ corresponds to the control system ${\mathbb L}$, using the notation in \S \ref{Lifting of singular trajectories.} applied to 
the projection $M = Z \to N = X$. 

In what follows, we take $Y$ an sufficiently small open neighbourhood of a point such that 
the composition $\pi_X\circ\pi_E$ of the bundle projection $\pi_E : E \to Z$ and $\pi_X : Z \to X$ is a fibration with a fibre 
which diffeomorphic to $\R^3$ (the product of $\R^2$ and an open interval). 

\bel
\label{immersion}
The differential map 
$\pi_{X*} \vert_{L \setminus 0} : L \setminus 0 \to TX$ 
restricted to $L$ off the zero section is an immersion. 
\enl

\Proof
We have two systems of local coordinates $y_1, \dots, y_5, z$ and $x_1, \dots, x_5, w$ of $Z$ such that 
$\pi_Y(y, z) = y$ and that $\pi_X(x, w) = x$. The line bundle $L$ is generated by $\zeta = \frac{\pa}{\pa z}$ 
and the line bundle $K$ is generated by $\varpi = \frac{\pa}{\pa w}$. 
By Lemma \ref{prolongation}, we have that $[\zeta, \varpi] \not= 0 \ \mod. E$. Therefore we have that 
$\ad(\zeta) : K \to TZ/K$ is injective. Therefore we have the result from Lemma \ref{immersion-general}. 
\QED


\ber 
{\rm
The linear hull $H(C) \subset TX$ of the cone structure $C \subset TX$ 
defines a contact structure $D'$ on $X$ (\cite{Bryant}). 
In fact the contact structure is defined by $D'_x = \pi_{X*}(E^{(4)}_z) (\subset T_xX)$, 
for any $z \in Z$ with $\pi_X(z) = x$ 
(cf. Lemma \ref{prolongation}). 
Moreover the tangent planes to the cones $C_{\gamma(t)}$, which is a Lagrangian plane in $D'_{\gamma(t)}$ 
along a singular trajectory $\gamma(t)$ of ${\mathbb C}$ 
is regarded as a curve in a Lagrange Grassmannian, which coincides with 
the (reduced) Jacobi curve introduced by Agrachev-Zelenko (cf.\cite{Zelenko2}). 
}
\enr

\bel
\label{lift}
Let $\gamma : [a, b] \to X$ be an immersive singular trajectory for the control system ${\mathbb C}$. 
Then there exists an immersive singular control $c : [a, b] \to L$ and 
a Lipschitz abnormal bi-extremal $\beta : [a, b] \to L\times_X T^*X$ for ${\mathbb C}$ 
corresponding to $\gamma$. 
\enl

\Proof
We keep to take the systems of coordinates used in Lemma \ref{lifting}, Lemma \ref{lifting2} 
and Lemma \ref{immersion}. 
The Hamiltonian $H : L\times_X T^*X \to \R$ is given, in this case, by 
$$
H(x, w; p; \lambda) = \lambda\langle p, \pi_{X*}(\zeta(x, w))\rangle, 
$$
where $(x, w, \lambda)$ (resp. $(x, w)$,  $\lambda$) is the system of coordinates of $L = \Ker(\pi_{Y*})$, (resp. of $Z$,  of $\pi_Y$-fibres) and $\zeta = \frac{\pa}{\pa z}$ as in Lemma \ref{immersion}. The control parameters are given by 
$w$ and $\lambda$. 
See the proofs of Lemma \ref{lifting} and Lemma \ref{lifting2}. 
Then the constrained Hamiltonian system is given by 
$$
\begin{array}{c}
\dot{x}(t) = \lambda(t)\pi_{X*}(\zeta(x(t), w(t))), 
\ \ 
\dot{p}_i(t) = - \lambda(t)\langle p(t), \frac{\pa}{\pa x_i}\pi_{X*}(\zeta(x(t), w(t)))\rangle, 
(1 \leq i \leq 5), 
\vspace{0.2truecm}
\\
\langle p(t), \pi_{X*}(\zeta(x(t), w(t)))\rangle = 0, \ \ 
\langle p(t), \frac{\pa}{\pa w}\pi_{X*}(\zeta(x(t), w(t)))\rangle = 0. 
\end{array}
$$
Note that $\pi_{X*}\circ\zeta : Z \to TX$ is a smooth vector field along the projection $\pi_X : Z \to X$. 
Take $w_0 \in \R$ such that $(x(a), w_0)$ belongs to the coordinate neighbourhood in $Z$ we are considering and 
$\dot{x}(a) \in (\pi_{X*})(L_{(x(a), w_0)})$. 
Since $\gamma(t) = x(t)$ is immersive, by Lemma \ref{immersion}, 
we have a unique smooth control $c : [a, b] \to L = \Ker(\pi_{Y*}) \subset E$, 
$c(t) = (x(t), w(t); \lambda(t))$ such that $\pi_L\circ c = \gamma$, $w(a) = w_0$. 
In fact, the coordinate $w(t)$ along $\pi_X$-fibre is determined by the velocity direction and the coordinate 
$\lambda(t)$ along $\pi_Y$-fibre is determined by the velocity of $\gamma$. 
Then $\beta(t) = (x(t), w(t); p(t);\lambda(t))$ is a Lipschitz abnormal bi-extremal. 
\QED

\

%

Thus we have reached to the stage to show the following: 

\bet
\label{duality}
A singular ${\mathbb C}$-path is the $\pi_X$-image of a $\pi_Y$-fibre. 
\ent

\Proof
Let $\gamma : I \to X$ be a singular ${\mathbb C}$-trajectory
Then, by Lemma \ref{lift}, $\gamma$ lifts to a Lipschitz 
abnormal bi-extremal $\beta : I \to L\times_X T^*X$ for ${\mathbb C}$. Then, by Lemma \ref{lifting}, 
$\gamma$ lifts to a singular ${\mathbb E}$-trajectory $\widetilde{\gamma}$ for the projection 
$\pi_X : Z \to X$, $\pi_X\circ \widetilde{\gamma} = \gamma$. 
Since $\gamma$ is immersive, $\widetilde{\gamma}$ is also immersive. 
By Proposition \ref{E-singular}, $\widetilde{\gamma}$ is either a $\pi_Y$-fibre or the
lift of a singular $D$-trajectory. In this case, $\widetilde{\gamma}$ must be a $\pi_Y$-fibre, 
and we have that $\gamma$ is the $\pi_X$-image of a $\pi_Y$-fibre up to parametrisation. 
\QED

\ber
{\rm 
For any point $x$ of $X$ and for any direction of $C_x$, there exists uniquely a singular ${\mathbb C}$-path through $x$ with the given direction. The space of singular ${\mathbb C}$-paths on $X$ is identified with $Y$. 
}
\enr

\ber
{\rm 
A singular ${\mathbb C}$-trajectory is not necessarily the image of a $\pi_Y$-fibre by $\pi_X$. 
For example a piecewise smooth ${\mathbb C}$-trajectory which consists of several images of $\pi_Y$-fibres by $\pi_X$ 
is also a singular ${\mathbb C}$-trajectory. This remark applies to singular ${\mathbb D}$-trajectories of course. 
}
\enr

\ber
{\rm 
The distribution $D \subset TY$ is also regarded as a cone structure defined by the control system 
${\mathbb K} : K = \Ker(\pi_{X*}) \subset TZ \to TY \to Y$. Then, naturally, singular paths for ${\mathbb K}$ 
coincide with singular $D$-paths by Lemmas \ref{lifting} and \ref{lifting2}. 
}
\enr

\ber
{\rm 
A singular bi-extremal for a control system is called {\it totally singular} if the Hessian matrix of the Hamiltonian 
with respect to the control parameter is zero along the bi-extremal and 
a trajectory is called {\it totally singular} if it possesses a totally singular bi-extremal (\cite{AS}). 
Then all singular trajectories on $X$ for ${\mathbb C}$ and on $Y$ for ${\mathbb D}$ are 
totally singular,  and therefore the duality of singular paths on $Y$ and $X$ 
holds on the level of totally singular paths as well. 
}
\enr

\section{Example: $G_2$-case.}

From the view point of twistor theory, 
the double fibration $Y \leftarrow Z \to X$ with $G_2$-symmetry
the $G_2$-Engel distribution $E \subset TZ$,
the Cartan distribution
$D \subset TY$ and the cone structure $C \subset TX$ have been explicitly constructed 
in \cite{IMT}. 
In this section, we give the explicit representations of the associated control systems
${\mathbb E} : E \hookrightarrow TZ \to Z$, ${\mathbb D} : D
\hookrightarrow TY \to Y$ and
${\mathbb C} : E \to TX \to X$, 
directly calculate the constrained Hamiltonian systems for the singular controls 
and determine the singular paths. 

\label{Example: $G_2$-case.}

\bee
{\rm 
(The control system ${\mathbb E}$ on $Z$.)
On a space $Z$ of dimension $6$ with coordinates $\lambda, x, y, z, u, v$, consider 
the explicit double fibrations
$$
(\lambda, x + \lambda y, y + \lambda z, v + \lambda x, u + \lambda(y^2 - xz)) 
\xleftarrow{\pi_Y} (\lambda, x, y, z, u, v) \xrightarrow{\pi_X} (x, y, z, u, v). 
$$
Then the $G_2$-Engel distribution $E = \Ker(\pi_{Y*})\oplus\Ker(\pi_{X*})$ is 
generated by 
$$
\xi_1 = \frac{\pa}{\pa \lambda}, \quad \xi_2 = \frac{\pa}{\pa z} - \lambda\frac{\pa}{\pa y} + 
\lambda^2\frac{\pa}{\pa x} - \lambda^3\frac{\pa}{\pa v} + 
(\lambda^3z + 2\lambda^2y + \lambda x)\frac{\pa}{\pa u}.
$$ 
For the coordinates $( \lambda, x, y, z, u, v;\kappa, p, q, r, \ell, m)$ of $T^{\ast}Z$, we have the Hamiltonian functions of $\xi_1, \xi_2$
$$
H_{\xi_1} = \kappa, \quad H_{\xi_2} = r - \lambda q + \lambda^2 p - 
\lambda^3 m + (\lambda^3 z + 2\lambda^2 y + \lambda x)\ell, 
$$
and the total Hamiltonian function
$
H = a H_{\xi_1} + b H_{\xi_2},
$
where $(a, b)$ is the fibre coordinate of $E$.
The constrained Hamiltonian system for singular controls is given by 
$$
\begin{array}{c}
\dot{\lambda} = a, \ \dot{x} = b\lambda^2, \ \dot{y} = - b\lambda, \ \dot{z} = b, \ 
\dot{u} = b(\lambda^3 z + 2\lambda^2 y + \lambda x), \ \dot{v} = - b\lambda^3, 
\\
\dot{\kappa} = - b\{-q+2\lambda p-3\lambda^2m+(3\lambda^2z+4\lambda y + x)\ell\}, 
\\
\dot{p} = -b\lambda\ell, \ \dot{q} = -b\cdot 2\lambda^2\ell, \ \dot{r} = -b\lambda^3\ell, 
\dot{\ell} = 0, \ \dot{m} = 0, 
\end{array}
$$
and $H_{\xi_1} = 0, H_{\xi_2} = 0$. 
We set $\ell = \ell_0, m = m_0$. 
By differentiating both sides of $H_{\xi_2}=0,$
we have 
$$
a\{-q+2\lambda p-3\lambda^2m_0+(3\lambda^2z+4\lambda y + x)\ell_0\} = 0.
$$

We consider two cases: the case where $a(t) \not= 0$ on a non-empty sub-interval and 
the case where $a(t)$ is identically zero. 

In the first case, we have
$-q+2\lambda p-3\lambda^2m_0+(3\lambda^2z+4\lambda y + x)\ell_0=0$ on the sub-interval. 
By differentiating both sides of the equation, 
we have 
$$
(2p - 6\lambda m_0) + (6\lambda z + 4y)\ell_0 = 0. 
$$
By differentiating both sides again,
we have 
$z\ell_0 - m_0 = 0$.
Assume $\ell_0 = 0$, then we have 
$p = q = r = 0$, which leads a contradiction with the condition on
abnormal bi-extremals.  Therefore $\ell_0 \not= 0$.
We have that $\dot{z} = 0$ therefore $b = 0,$ 
and that $x, y, z, u, v$ are constant on the sub-interval. 
Moreover, if the singular trajectory is immersive, we have that $x, y, z, u, v$ are constant on the whole interval, similarly to the proof of Proposition \ref{E-singular}. 
Thus, in this case the singular control gives a $\pi_X$-fibre. 

If $a$ is identically zero, then $\lambda$ is constant and 
$$
(x + \lambda y)' = 0, (y + \lambda z)' = 0, (v + \lambda x)' = 0, (u + \lambda(y^2 - xz))' = 0. 
$$
Therefore, in the second case, the singular control gives a $\pi_Y$-fibre. 

}
\ene

\bee
{\rm 
(The control system ${\mathbb D}$ on $Y$.)
For the Cartan distribution $D$ on $Y$, 
let us consider the space with coordinates $\lambda, \mu, \nu, \tau, \sigma$ and 
the distribution $D \subset TY$ generated by
$$
\eta_1 = \frac{\pa}{\pa \lambda} + \nu\frac{\pa}{\pa \mu} - (\lambda\nu - \mu)\frac{\pa}{\pa \tau} 
+ \nu^2\frac{\pa}{\pa \sigma}, \quad 
\eta_2 = \frac{\pa}{\pa \nu} - \lambda\frac{\pa}{\pa \mu} + \lambda^2\frac{\pa}{\pa \tau} - 
(\lambda\nu + \mu)\frac{\pa}{\pa \sigma}. 
$$
For the coordinates $( \lambda, \mu, \nu, \tau, \sigma;p, q, r, \ell, m)$ of $T^*Z$, we set 
$$
H_{\eta_1} = p + \nu q - (\lambda\nu - \mu)\ell + \nu^2m, \quad 
H_{\eta_2} = r - \lambda q + \lambda^2 \ell - (\lambda\nu + \mu)m, 
$$
and $H = u_1H_{\eta_1} + u_2H_{\eta_2}$,
where $(u_1,u_2)$ is the fibre coordinate of $D$. Then we have the constrained Hamiltonian equation: 
$$
\begin{array}{c}
\dot{\lambda} = u_1, \ \dot{\mu} = u_1\nu - u_2\lambda, \dot{\nu} = u_2, \ 
\dot{\tau} = -u_1(\lambda\nu - \mu) + u_2\lambda^2, 
\ 
\dot{\sigma} = u_1\nu^2 - u_2(\lambda\nu + \mu), 
\\
\dot{p} = u_1\nu\ell - u_2(2\lambda\ell - \nu m), \ 
\dot{q} = - u_1\ell + u_2 m, 
\ 
\dot{r} = - u_1(q - \lambda\ell + 2\nu m) + u_2\lambda m, 
\\
\dot{\ell} = 0, \ \dot{m} = 0, \quad 
p + \nu q - (\lambda\nu - \mu)\ell + \nu^2m = 0, \ 
 r - \lambda q + \lambda^2 \ell - (\lambda\nu + \mu)m = 0. 
\end{array}
$$
Then by direct calculations we can conclude that singular $D$-paths are $\pi_Y$-projections of $\pi_X$-fibres. 
Because this claim is known, we omit the details. 
}
\ene

\bee
{\rm 
(The control system ${\mathbb C}$ on $X$.) 
Let $\lambda,x,y,z,u,v$ be local coordinates of  $Z$.
From the explicit form of the projection $\pi_Y$, we see that the line
bundle $L = \Ker(\pi_{Y*})$ is generated by
$$
  \frac{\pa}{\pa z} + \lambda(- \frac{\pa}{\pa y} + x \frac{\pa}{\pa u})
  + \lambda^2(\frac{\pa}{\pa x} + 2y \frac{\pa}{\pa u}) + \lambda^3(-
 \frac{\pa}{\pa v} + z \frac{\pa}{\pa u}).
$$
Then we can regard the cone structure $C\subset TX$
as the control system ${\mathbb C} : L \to TX \to X$ on $X$ 
given by
$$
 \dot{\gamma} = \mu\left\{ \frac{\pa}{\pa z} + \lambda(- \frac{\pa}{\pa y}
 + x \frac{\pa}{\pa u})
 + \lambda^2(\frac{\pa}{\pa x} + 2y \frac{\pa}{\pa u}) + \lambda^3(-
 \frac{\pa}{\pa v} + z \frac{\pa}{\pa u})\right\},
 $$
where $\gamma(t) = (x(t), y(t), z(t), u(t), v(t))$ and
$\lambda, \mu$ are regarded as control parameters.
The Hamiltonian function is given by 
$$
H = \mu\{ r + \lambda(- q + x\ell) + \lambda^2(p + 2y\ell) + \lambda^3(-m + z\ell)\},
$$
where $(x, y, z, u, v; p, q, r, \ell, m)$ is a coordinate system of $T^*X$. 
Then the constrained Hamiltonian system is given by 
$$
\begin{array}{c}
\dot{x} = \mu\lambda^2, \ \dot{y} = - \mu\lambda, \ \dot{z} = \mu, \ 
\dot{u} = \mu(\lambda x + 2\lambda^2 y + \lambda^3 z), \ \dot{v} = - \mu\lambda^3, 
\\
\dot{p} = - \mu\lambda\ell, \ \dot{q} = - 2\mu\lambda^2\ell, \ \dot{r} = - \mu\lambda^3\ell, \ \dot{\ell} = 0, \ \dot{m} = 0, 
\\
r + \lambda(- q + x\ell) + \lambda^2(p + 2y\ell) + \lambda^3(-m + z\ell) = 0, 
\\
\mu(-q + x\ell + 2\lambda(p + 2y\ell) + 3\lambda^2(-m + z\ell) = 0. 
\end{array}
$$
Set $\ell = \ell_0, m = m_0$. 
Since our trajectory is immersive, we have $\mu \not= 0$. 
Therefore we have 
$$
-q + x\ell_0 + 2\lambda(p + 2y\ell_0) + 3\lambda^2(-m_0 + z\ell_0) = 0
$$
By differentiating both sides, we have 
$$
\dot{\lambda}\{ 2(p + 2y\ell_0) + 6\lambda(-m_0 + z\ell_0)\} = 0. 
$$
We assume $\dot{\lambda} \not= 0$. Then 
$$
2(p + 2y\ell_0) + 6\lambda(-m_0 + z\ell_0) = 0. 
$$
By differentiating both sides, we have 
$\dot{\lambda}(- m_0 + z\ell_0) = 0$, so $- m_0 + z\ell_0 = 0$. By differentiating both sides, 
we have $\mu\ell_0 = 0$, so $\ell = \ell_0 = 0$. Then we have 
$m = m_0 = 0$, and then $p = q = r = 0$. This contradicts with the condition on abnormal bi-extremals. 
Therefore we have $\dot{\lambda} = 0$, and $\lambda$ is constant. Set $\lambda = \lambda_0$. 
Then we have, up to parametrisation and orientation, 
$$
x(t) = \lambda_0^2 t + x_0, \ y(t) = - \lambda_0  t + y_0, \ z(t)
=t + z_0, \ v(t) = -\lambda_0^3 t + v_0, 
$$
and 
$$
u(t) = (\lambda_0 x_0 + 2\lambda_0^2 y_0 + \lambda_0^3 z_0)t+u_0, 
$$
where $x_0\,\,(\mbox{resp.}\,\,y_0,z_0,v_0,u_0)$ is the initial value of $x(t)\,\,(\mbox{resp.}\,\,y(t),z(t),v(t),u(t))$.
Thus we have shown that any singular ${\mathbb C}$-path is the $\pi_X$-projection of a $\pi_Y$-fibre in $G_2$-case 
directly. 
}
\ene

\ber
{\rm 
In $G_2$-case \cite{IMT}, 
the $\pi_X$-image of a $\pi_Y$-fibre, namely a singular ${\mathbb C}$-path in $X$, 
is called a {\it Monge line}, while 
the $\pi_Y$-image of a $\pi_X$-fibre, namely a singular ${\mathbb D}$-path in $Y$, 
is called a {\it Cartan line}. 
}
\enr

\

\begin{flushleft}
Goo ISHIKAWA, \\
Department of Mathematics, Hokkaido University, 
Sapporo 060-0810, Japan. \\
e-mail : ishikawa@math.sci.hokudai.ac.jp \\

\

Yumiko KITAGAWA, \\
Oita National College of Technology, 
Oita 870-0152, Japan. \\
e-mail : kitagawa@oita-ct.ac.jp

\

Wataru YUKUNO, \\
Department of Mathematics, Hokkaido University, 
Sapporo 060-0810, Japan. \\
e-mail : yukuwata@math.sci.hokudai.ac.jp

\end{flushleft}

\end{document}